# Finite volume methods for elasticity with weak symmetry

Eirik Keilegavlen[1] and Jan Martin Nordbotten[1,2]
1. Department of Mathematics, University of Bergen
2. Department of Civil and Environmental Engineering, Princeton University

Corresponding author: Eirik Keilegavlen, email: Eirik.Keilegavlen@math.uib.no

## Abstract
We introduce a new cell-centered finite volume discretization for elasticity with weakly enforced symmetry of the stress tensor. The method is motivated by the need for robust discretization methods for deformation and flow in porous media, and falls in the category of multi-point stress approximations (MPSA). By enforcing symmetry weakly, the resulting method has additional flexibility beyond previous MPSA methods. This allows for a construction of a method which is applicable to all grid types, and in particular the method amends a crucial shortcoming in previous MPSA methods for simplex grids.

By formulating the method as a discrete variational problem, we prove convergence of the new method for a wide range of problems, with conditions that can be verified at the time of discretization. We present the first set of comprehensive numerical tests for the MPSA methods in three dimensions, covering Cartesian and simplex grids, with both heterogeneous and nearly incompressible media. The tests show that the new method consistently is second order convergent in displacement, despite being lowest order, with a rate that mostly is between 1 and 2 for stresses. The results further show that the new method is more robust and computationally cheaper than previous MPSA methods.

## 1. Introduction

Fluid flow in porous materials is intrinsically coupled to mechanical stresses induced on the skeleton. In many geological subsurface applications, such as $CO_2$ storage [1], geothermal energy extraction, e.g. [2], and petroleum extraction [3], engineering design calls for high flow rates, and consequently the induced mechanical stresses become significant. In this context, simulation of coupled flow and deformation in geological porous media is becoming increasingly important.

Numerical simulation of the Biot equations is challenging due to the different mathematical nature of the flow and mechanical equations. In order to allow stable discretizations in the presence of multiple fluid phases, the flow equations are naturally discretized by conservative schemes, such as finite volume methods [4]. In contrast, it is customary to discretize the mechanical equations using finite element methods, e.g. [5]. This situation has the disadvantage that finite volume and finite element methods inherently use different data structures and are best adapted to different grid types. This makes the construction of efficient simulation codes difficult, and as a consequence fixed point iteration between disparate computational tools has become industry standard [6].

The status quo for numerical simulation of the Biot equations is not optimal, and several authors have tried to address this aspect. Thus, simulation approaches built on finite elements [7], finite differences [8] and finite volumes [6], [9] have recently been considered. We argue that in the context of applications related to multi-phase flow in geological porous media, it is imperative to use conservative discretizations for the flow equations [10], [11]. As such, we see a need for finite volume methods for mechanics which are designed to handle the grids and material discontinuities typical of industrial reservoir simulation [12]–[14].

Recently, we proposed a family of finite volume method for mechanics, referred to as Multi-Point Stress Approximation (MPSA) methods [15]. These methods build on the same concepts and the same data structures as the Multi-Point Flux Approximation (MPFA) methods common for multi-phase flows in porous media [6], [16], [17]. Since introducing the MPSA methods, we have demonstrated their applicability to the coupled Biot system [18].

We have recently shown that the MPSA-O discretizations convergent discretizations for mechanics [19]. The proof highlights the role of local coercivity conditions, which are a functions of the local geometry, material parameters, and discretization scheme. These local coercivity conditions can be verified for many classes of grids relevant for applications (Cartesians and so-called PEBI or Voronoi grids). However, the situation for simplex grids is less clear, and the existing MPSA-type discretizations may fail unless strict conditions of acute angles are enforced on the grids, even in 2D. Thus the situation is that the MPSA methods until now have been more restrictive in terms of grids than their flow counterpart MPFA.

The existing MPSA methods are all generalizations of their MPFA counterparts. Herein, we present the first MPSA method which is explicitly defined to exploit the particular nature of the equations of elasticity. In particular, it is well known in the context of mixed finite elements for elasticity that obtaining inf-sup stable discretizations may require quite sophisticated spaces [20]. On the other hand, imposing symmetry in an integral sense may lead to simpler discretization methods, as studied in the setting of both mixed finite elements [21] and a mimetic finite differences [22]. This concept is also referred to as weak symmetry, which is the key tool we will use to obtain the discretization proposed herein.

We are thus motivated by this observation to seek a finite volume discretization with weakly enforced symmetry, particularly tailored to discretization of elasticity. We refer to this new method as MPSA-W, to emphasize that the symmetry is in a weak sense. In exchange for the reduced symmetry, several practical and computational advantages are realized. In particular, the MPSA-W method improves over existing MPSA methods as it i) avoids solving computationally expensive local minimization problems at the time of discretization, ii) is robust for simplex grids. It shares with existing MPSA discretizations the applicability to general polygonal grids, material discontinuities and Poisson ratio. Furthermore, despite the nomenclature, all MPSA methods, including the MPSA-W method proposed herein, allow for a local post-processing which recovers a point-wise symmetric stress tensor.

We structure the remaining presentation as follows. In Section 2, we first present the governing equations, and then introduce the notation necessary to describe the discretization method. The MPSA-W method is presented in two different settings: Section 3 describes a framework which is natural for implementation, and also shows the similarities with the related method for the scalar equation. In Section 4, the method is cast in a discrete variational framework which allows us to discuss the theoretical convergence properties of the method. In Section 5, the performance of the MPSA-W method is assessed by numerical convergence tests in 2D and 3D, before Section 6 provides concluding remarks.

## 2. Equations and grid

In this section we introduce the governing equations and notation necessary to describe the discretization schemes.

### 2.1 Governing equations

We consider the static momentum balance equation for an elastic medium, which in Lagrangian coordinates reads

$$\int_{\partial \Omega} \boldsymbol{T}(\boldsymbol{n}) dA + \int_\Omega \boldsymbol{f} dV = 0 \qquad (2.1.1)$$

Where $\boldsymbol{T}(\boldsymbol{n})$ are the surface traction vectors on the boundary of some domain $\Omega \in R^D, D = \{1,2,3\}$, with normal vector $\boldsymbol{n}$ and $\boldsymbol{f}$ represents body forces. Assuming small deformations and a linear stress-strain relation, the surface forces can be expressed as

$$\boldsymbol{T}(\boldsymbol{n}) = \boldsymbol{\pi} \cdot \boldsymbol{n} = (\mathbb{C} : \boldsymbol{\epsilon}) \cdot \boldsymbol{n} \qquad (2.1.2)$$

Here, $\boldsymbol{\pi}$ is the 1st Piola-Kirchhoff stress tensor, $\mathbb{C}$ is the stiffness tensor, and $\boldsymbol{\epsilon}$ is the symmetric part of the deformation gradient $\boldsymbol{\epsilon} = (\nabla \boldsymbol{u} + (\nabla \boldsymbol{u})^T)/2$. Since we are in the context of small deformations in Eucledian space, the Piola-Kirchhoff stress tensors coincide with the Cauchy stress tensor, and are all symmetric.

In the case of an isotropic media, $\mathbb{C}$ can be expressed by the Lamé parameters $\mu, \lambda$ to arrive at the familiar form

$$\boldsymbol{\pi} = 2\mu\boldsymbol{\epsilon} + \lambda\, tr(\boldsymbol{\epsilon})\boldsymbol{I} \qquad (2.1.3)$$

When the solution has sufficient smoothness, Equations (2.1.1)-(2.1.2) can be re-written as the differential equation

$$\nabla \cdot \boldsymbol{\pi} = \boldsymbol{f} \qquad (2.1.4)$$

Thus Equations (2.1.1)-(2.1.3) are equivalent with Equations (2.1.3)-(2.1.4).

### 2.2 Grid

To define our numerical method for solving (2.1.1)-(2.1.3), we need notation to describe the grid and associated quantities. The notation introduced below is identical to the one used in [19], as this will allow us to adapt convergence results therein to the weakly symmetric discretization.

We denote a generic cell $K$, with volume (area in 2D) $m_K$, and let $\mathcal{T}$ be the set of all cells. Similarly, $\sigma$ denotes a generic face, with area $m_\sigma$, and $\mathcal{F}$ is the set of all faces in the grid. For $D = 3$, we will further denote a generic edge in the grid $E$. Finally, $\mathcal{V}$ is the set of vertexes in the grid, and $s$ will denote a generic vertex.

To simplify the description of the method, we also introduce notation to represent the relation between cells, faces, edges and vertexes. The faces on the boundary of a cell $K$ is denoted $F_K$, so that $\partial K = \bigcup_{\sigma \in \mathcal{F}_K} \sigma$, and the vertexes of $K$ is $\mathcal{V}_K$. Similarly, for each face, the neighboring cells are denoted $\mathcal{T}_\sigma$ and the face corners for $\mathcal{V}_\sigma$. For each vertex $s \in \mathcal{V}$, we denote the adjacent cells by $\mathcal{T}_s$ and the adjacent faces by $\mathcal{F}_s$.

We introduce the following geometric quantities: Let $\boldsymbol{x}_K$ and $\boldsymbol{x}_\sigma$ represent cell and face centers, respectively. For each face, we let $\boldsymbol{n}_{K,\sigma}$ denote the unit normal vector which is outer with respect to the cell $K \in \mathcal{T}_\sigma$, and let $d_{K,\sigma}$ be the distance from $\boldsymbol{x}_K$ to $\boldsymbol{x}_\sigma$. Finally, we denote by $\boldsymbol{x}_s$ vertex coordinates.

Our finite volume method will be constructed to have conservation properties on the cells $\mathcal{T}$. However, much of the construction of the discretization is based on a further partitioning of the cells into sub-cells, see Figure 1. With each vertex $s \in \mathcal{V}_K$ we associate a sub-cell which we identify by the pair $(s, K)$, with a volume $m_K^s$, so that $\sum_{s \in \mathcal{V}_K} m_K^s = m_K$. Similarly, sub-faces are defined for each vertex $s \in \mathcal{V}_\sigma$, identified by the pair $(s, \sigma)$, and with an area $m_\sigma^s$.

With the above notation, momentum conservation stated for $K$ reads

$$-\int_K \boldsymbol{f} dV = \int_{\partial K} \boldsymbol{T}(\boldsymbol{n}) dA = \sum_{\sigma \in \mathcal{F}_K} \boldsymbol{T}_{K,\sigma} \tag{2.2.1}$$

where $\boldsymbol{T}_{K,\sigma}$ denotes the forces acting on $\sigma$.

## 3. Discretization

This section is devoted to the introduction of MPSA methods with both weak and strong symmetry. The presentation and notation in this section is chosen so that it is convenient for implementation, while Section 4 offers an alternative approach which is better suited for convergence analysis.

### 3.1 Discrete variables

We discretize the displacements in terms of cell center displacements, and denote these degrees of freedom $\boldsymbol{u}_K$. To each cell $K$, and for each component of the displacement vector, we further associate basis functions, which are piece-wise linear functions on each sub-cell, and with support on all sub-cells for which $\boldsymbol{x}_K$ is a vertex, see Figure 1a. The gradients that define the basis functions of cell $K$ are denoted $\boldsymbol{G}_K$, and together with $\boldsymbol{u}_K$ they define the discrete displacements everywhere. To obtain a global discretization in terms of cell center displacements, the gradients must be eliminated by local calculations. Before we carry out the elimination, it is therefore necessary to introduce the discretized tractions on the surface, and in particular introduce the notion of weakly enforced symmetry of the stress tensor.

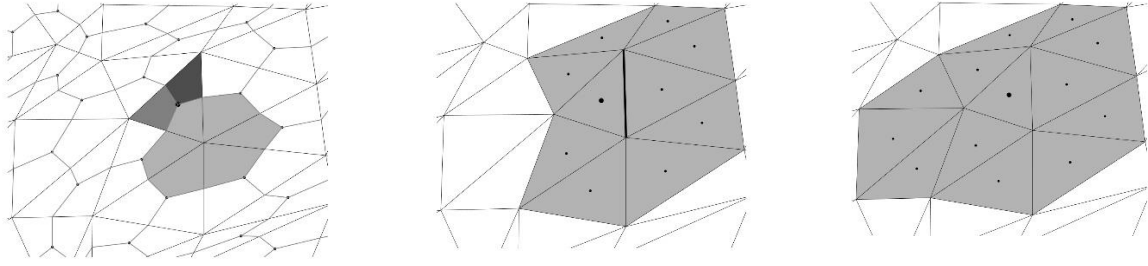

| a) Sub-cells | b) Stress stencil | c) Displacement stencil |

**Figure 1:** Grid structure used in the discretization of the center cell (marked with bullet). Left: Partitioning of triangular cells into sub-cells; Sub-cells with light gray share the same vertex, and form the domain of a local calculation. Middle: Stencil for stress discretization on one face; Right: Stencil for displacement discretization for the whole cell.

## 3.2 Discrete representation of face stresses

In the continuous setting, the stress at a point is given by $\boldsymbol{\pi} = \mathbb{C}_K : \boldsymbol{\epsilon}(\boldsymbol{u})$. The discrete analogue can be written in terms of the discrete gradients $\boldsymbol{G}_K^s$ as

$$\underline{\boldsymbol{\pi}}_{K,s} = \frac{\mathbb{C}_K : \boldsymbol{G}_K^s + (\mathbb{C}_K : \boldsymbol{G}_K^s)^T}{2} \tag{3.2.1}$$

Where we have exploited that $\mathbb{C}_K$ can written as a symmetric tensor. This stress expression (3.2.1) was introduced in [15], and leads to the MPSA-O methods with strong symmetry (emphasized by the underscore).

The purpose of this paper is to explore an alternative formulation of stresses, which leads to a method with weak symmetry. We will see that this method has superior qualities. To this end, the discretized term $(\mathbb{C}_K : \boldsymbol{G}_K^s)^T$ is replaced by a quantity averaged over surrounding sub-cells. Specifically, we introduce two notions of weak symmetry, one associated with vertexes of the grid, while the second is associated with edges of the grid. These notions coincide in 2D, but lead to distinct methods in 3D.

We introduce the notation for the average associated with vertexes, $\langle \cdot \rangle_s$, and introduce the discrete stress

$$\boldsymbol{\pi}_{K,s} = \frac{\mathbb{C}_K : \boldsymbol{G}_K^s + \langle \mathbb{C}_K : \boldsymbol{G}_K^s \rangle_s^T}{2} \tag{3.2.2}$$

We note that the term weak symmetry is justified since

$$2\langle \boldsymbol{\pi}_{K,s} - \boldsymbol{\pi}_{K,s}^T \rangle_s = \langle \mathbb{C}_K : \boldsymbol{G}_K^s + \langle \mathbb{C}_K : \boldsymbol{G}_K^s \rangle_s^T - (\mathbb{C}_K : \boldsymbol{G}_K^s)^T - \langle \mathbb{C}_K : \boldsymbol{G}_K^s \rangle_s \rangle_s = 0$$

by linearity of the averaging operator and since $\langle \langle \mathbb{C}_K : \boldsymbol{G}_K^s \rangle_s \rangle_s = \langle \mathbb{C}_K : \boldsymbol{G}_K^s \rangle_s$.

While edges and vertexes coincide in 2D, they are distinct in 3D. A slightly stronger notion of symmetry can thus be enforced at the expense of the additional data-structure needed to handle edges. We introduce the notation for edge average $\langle \cdot \rangle_e$, and introduce the discrete stress

$$\boldsymbol{\pi}_{K,s,e} = \frac{\mathbb{C}_K : \boldsymbol{G}_K^s + \langle \mathbb{C}_K : \boldsymbol{G}_K^s \rangle_e^T}{2} \tag{3.2.3}$$

Again the term weak symmetry is justified as above. We will refer to the stress discretization (3.2.2) as the MPSA-W method with vertex-wise averaging, while (3.2.3) will be termed edge-wise averaging; in 2D we will simply refer to the MPSA-W method.

Based on any of the discrete stresses $\chi_{K,s,e} = \{\underline{\pi}_{K,s}, \pi_{K,s}, \pi_{K,s,e}\}$, we can define discrete tractions according to:

$$T^{\sigma}_{K,s} = \sum_{e \in \mathcal{E}_s \cap \mathcal{E}_\sigma} m^{s,e}_\sigma \chi_{K,s,e} \cdot n_{K,\sigma} \tag{3.2.4}$$

## 3.3 Local linear system

Independent of which discretization of stress is chosen, the finite volume framework requires that the stress over a sub-face is continuous, that is

$$T^{\sigma}_{R,s} = -T^{\sigma}_{L,s}, \qquad \text{with } \{R, L\} = \mathcal{T}_\sigma \tag{3.3.1}$$

We note that via the definition of discrete tractions (3.2.1)-(3.2.3), this is an equation for gradients $G^s_K$. For implementation purposes, it is important to note that by construction the term

$$\langle \mathbb{C} : G^S_K \rangle^T_s \cdot n_{K,\sigma} = -\langle \mathbb{C} : G^S_L \rangle^T_s \cdot n_{L,\sigma} \tag{3.3.2}$$

(and similarly for the average over $e$), thus only the part $\mathbb{C} : G^S_K$ will contribute to (3.3.1). This is significant, since this gives a locally coercive system, which significantly simplifies the method.

We will also require certain continuity properties for the displacement: For each sub-face $(s, \sigma)$, we require the displacement to be continuous in a point $x^s_{\sigma,C} = \eta x_s + (1-\eta) x_\sigma$. Here, $\eta$ is a parameter, which based on experience from MPFA-O methods we will set to 1/3 on simplex grids and 0 otherwise [23], [24], see comments towards the end of Section 4. For each face, this gives an equation

$$u_R + G^S_R \left( x^s_{\sigma,C} - x_R \right) = u_L + G^S_L \left( x^s_{\sigma,C} - x_L \right), \text{ with } \{R, L\} = \mathcal{T}_\sigma \tag{3.3.3}$$

Note that the use of a single continuity point leads to the MPSA-O method when using the discrete stress $\underline{\pi}_{K,s}$. However, this method is only applicable to simplex grids, thus the *generalized* MPSA-O method (with a minimization problem for discontinuities) is used for general grids [15]. The MPSA-W method avoids this additional complexity, and we will limit the presentation to this simpler case.

While (3.3.3) enforces a certain displacement continuity between subcells located in different cells, we require no such continuity between subcells contained within a single cell $K$. This is a key component of the MPSA methods, which allows us to write the gradients $G^S_K$ in terms of cell-center displacements in Equation (3.2.4) by a local calculation, involving only subcells surrounding a single vertex. To be specific, the stress continuity (3.3.1) and the displacement continuity (3.3.3) can be combined into a linear system

$$\begin{pmatrix} S_G & 0 \\ D_G & I \\ 0 & I \end{pmatrix} \begin{pmatrix} G \\ U \end{pmatrix} = \begin{pmatrix} 0 \\ 0 \\ I \end{pmatrix} \tag{3.3.4}$$

Here, $S_G$ represents products on the from $n^T \mathbb{C}$, $D_G$ contains distances from cell centers to continuity points, while $I$ is the identify matrix. The unknowns are all gradients in the subcells ($G$) and all cell center displacements ($U$). The first two rows in (3.3.4) thus represent continuity of forces and displacement, while each row in the lower block enforces a unit displacement in one cell-center at a time. When (3.3.4) is solved for $G$, and the gradients are inserted into (3.2.4), we obtain the desired expressions for surface forces on a sub-face in terms of cell center displacements. The unique solvability of Equations (3.3.4) is an important aspect of the method. We return to this aspect in Section 4.

We note that boundary conditions are easily incorporated into (3.3.4) by adjusting either (3.3.1) (for Neumann conditions), (3.3.3) (Dirichlet conditions) or a combination thereof (mixed conditions).

## 3.4 Global linear system

The local calculation is carried out for all faces in the grid to obtain discrete traction in terms of the surrounding displacements, which we will write on the form $T^\sigma_{K,s} = \sum_{K' \in \mathcal{T}_s} t_{K,K',s,\sigma} \, \boldsymbol{u}_{K'}$. Insertion into (2.1.2) now gives

$$-\int_K \boldsymbol{f} \, dV = \sum_{\sigma \in \mathcal{F}_K} \boldsymbol{T}_{K,\sigma} = \sum_{\sigma \in \mathcal{F}_K} \sum_{s \in \mathcal{V}_\sigma} \boldsymbol{T}^\sigma_{K,s} = \sum_{\sigma \in \mathcal{F}_K} \sum_{s \in \mathcal{V}_\sigma} \sum_{K' \in \mathcal{T}_s} t_{K,K',s,\sigma} \, \boldsymbol{u}_{K'}$$

(3.4.1)

This expresses force balance on cell $K$ in terms of the displacements of all cells sharing at least one vertex with $K$. When (3.4.1) is assembled for all $K \in \mathcal{T}$, and boundary conditions are taken into account, this results in a linear system for the cell center displacements.

## 4. Convergence

Convergence of the generalized MPSA-O method has been proved, conditional on locally computable criteria relating to grid and material parameters [19]. We will cast the MPSA-W method in the same framework. However, we will see that the locally criteria are different, and indeed favor the MPSA-W method over the MPSA-O methods.

## 4.1 Discrete variational framework

Equations (3.2.4) and (3.4.1) are equivalent to a discrete variational form. Indeed, let us consider the function space $\mathcal{H}_\mathcal{T}$ consisting of cell-center displacements, and the space $\mathcal{H}_\mathcal{C}$ consisting of cell-center displacements as well as displacements at the face centers [19], [25], [26]. Then on the space $\mathcal{H}_\mathcal{C}$ we can define a gradient which is consistent with a planar interpolation between the cell-center variables and the face-center variables:

$$(\overline{\nabla} \boldsymbol{u})^s_K = \sum_{\sigma \in \mathcal{F}_K \cap \mathcal{F}_s} (\langle \boldsymbol{u} \rangle^\sigma_{K,s} - \boldsymbol{u}_K) \otimes \boldsymbol{g}^s_{K,\sigma} \qquad (4.1.1)$$

The vectors $\boldsymbol{g}^s_{K,\sigma}$ are uniquely defined in 2D and for grids where for each cell, no more than three of its faces meet in a vertex in 3D [26] by the equations

$$I = (\overline{\nabla}x)_K^S = \sum_{\sigma \in \mathcal{F}_K \cap \mathcal{F}_s}(\langle x \rangle_{K,s}^{\sigma} - x_K) \otimes g_{K,\sigma}^S \tag{4.1.2}$$

We see from Equation (3.3.3) that this gradient is consistent with $G_L^S$. In order to enforce a conservation principle, we additionally utilize a gradient built from the normal vectors of the grid, and which corresponds to a divergence in the variational framework:

$$(\widetilde{\nabla}u)_K^S = \frac{1}{m_K^S}\sum_{\sigma \in \mathcal{F}_K \cap \mathcal{F}_s} m_\sigma^s(\langle u \rangle_S^\sigma - u_K) \otimes n_{K,\sigma} \tag{4.1.3}$$

With these discrete gradient operators, we can re-write the finite volume method for elasticity as a variational problem by introducing the bilinear form $(u, v) \in \mathcal{H}_C \times \mathcal{H}_C$

$$b_\mathcal{D}(u, v) = \sum_{K \in \mathcal{T}} \sum_{s \in \mathcal{V}_K} m_K^S [\mathbb{C}_K : (\overline{\nabla}u)_K^S + \langle \mathbb{C}_K : (\overline{\nabla}u)_K^S \rangle_x^T] : (\widetilde{\nabla}v)_K^S \tag{4.1.4}$$

Note that $\langle \mathbb{C}_K : (\overline{\nabla}u)_K^S \rangle_x$ corresponds to the average stress around vertex or edge, with $x \in \{s, e\}$. This corresponds to enforcing the symmetry of the stress tensor weakly. The equivalent variational problem to equations (3.2.4) and (3.4.1) is then stated as: Find $u_\mathcal{D} \in \mathcal{H}_C$ such that

$$b_\mathcal{D}(u_\mathcal{D}, v) = \int_\Omega f \cdot \mathcal{P}_{C,\mathcal{T}} v \, dx \qquad \text{for all } v \in \mathcal{H}_C \tag{4.1.5}$$

Here, $\mathcal{P}_{C,\mathcal{T}}$ is the projection of $\mathcal{H}_C$ onto piecewise constants on the cells. Equivalence can easily be verified by considering the canonical basis for the test functions $v$ [19], [26]. It is particularly useful to split Equation (4.1.5) into two components, in the sense of a variational multi-scale method [27], associated with the cell-center space $\mathcal{H}_\mathcal{T}$ and the cell-face space $\mathcal{H}_\mathcal{F} = \mathcal{H}_C/\mathcal{H}_\mathcal{T}$: Find $(u_\mathcal{T}, u_\mathcal{F}(u_\mathcal{T})) \in \mathcal{H}_\mathcal{T} \times \mathcal{H}_\mathcal{F}$ such that

$$b_\mathcal{D}(u_\mathcal{T} + u_\mathcal{F}(u_\mathcal{T}), v) = \int_\Omega f \cdot \mathcal{P}_{C,\mathcal{T}} v \, dx \qquad \text{for all } v \in \mathcal{H}_\mathcal{T} \tag{4.1.6}$$

$$b_\mathcal{D}(u_\mathcal{F}(u_\mathcal{T}), v) = -b_\mathcal{D}(u_\mathcal{T}, v) \qquad \text{for all } v \in \mathcal{H}_\mathcal{F} \tag{4.1.7}$$

The point is as noted in Section 3.3, that Equation (4.1.7), and thus the calculation of the projection operators $\Pi_{FV} u_\mathcal{T} = u_\mathcal{T} + u_\mathcal{F}(u_\mathcal{T})$, can be resolved purely locally around each vertex of the grid. Thus Equation (4.1.6) represents the finite volume method with degrees of freedom in the cell centers only.

In order to invoke the convergence results of [19], we need that the local systems used to calculate $\Pi_{FV}$ have a unique solution. In general, this appears to always hold true in practice [26], [28], however the proof is elusive even in the case of the scalar MPFA-O method. Here, we thus simply that the local problems for the MPSA-W method are essentially equivalent to the MPFA-O method.

**Lemma 4.1** Let $\lambda = 0$. Then the local problems (4.1.7) for the MPSA-W methods are equivalent to those of the MPFA-O method.

*Proof.* Due to (3.3.2), the local problems for the case of $\lambda = 0$ are nothing but three independent scalar problems, each of which is equivalent to the MPFA-O methods. ∎

It is important to note that this relationship does not hold for the MPSA-O method, since we do not have (3.3.2). This necessitated the development of the generalized MPSA-O method [15], for which Lemma 4.1 can be shown due to the existence of a minimization problem. Thus Lemma 4.1 is one of the main advantages of the method proposed herein.

Furthermore, we need to verify that the bilinear forms are sufficiently close to being symmetric within the subspace defined by the solution to the local problems. We denote by $b_{\mathcal{D},s}$ and $\Pi_{FV,s}$ the restriction of the bilinear form and the projection operator to the cells $K \in \mathcal{T}_s$:

**Local Coercivity Condition:** For every vertex $s \in \mathcal{V}$, there exists a constant $\theta_{2,s} \geq \theta_2 > 0$ such that the bilinear form $b_{\mathcal{D},s}$ and the projection $\Pi_{FV,s}$ satisfy for all $\boldsymbol{u} \in \Pi_{FV,s}\mathcal{H}_{\mathcal{T}}/\Re(\Omega)$

$$b_{\mathcal{D},s}(\boldsymbol{u},\boldsymbol{u}) \geq \theta_{2,s}|\boldsymbol{u}|^2_{b_{\mathcal{D},s}}$$

Where the local energy semi-norm is associated with the symmetrized bilinear form

$$|\boldsymbol{u}|^2_{b_{\mathcal{D},s}} = \sum_{K \in \mathcal{T}_s} m_K^s \left(\overline{\underline{\nabla}}\boldsymbol{u}\right)^s_K : \mathbb{C} : \left(\overline{\underline{\nabla}}\boldsymbol{u}\right)^s_K$$

This assumption can be verified locally while assembling the discretization, and moreover it can be verified *a priori* for certain classes of meshes [19]. We will return to both the solvability of the local system as well as the local coercivity condition the results section.

Given the solvability of the local problems and the local coercivity condition, convergence of the method follows directly from the analysis in [19]. We recall the three main results:

**Theorem 4.1:** For given parameter field $\mathbb{C}$, and mesh $\mathcal{D}$, let the local coercivity condition hold. Then the coarse variational problem (4.1.6) is coercive in the sense that it satisfies if a) $\Gamma_D$ is measurable then for all $\boldsymbol{u} \in \mathcal{H}_{\mathcal{T}}$ and b) if $\Gamma_N = \partial\Omega$ then for all $\boldsymbol{u} \in \mathcal{H}_{\mathcal{T}}/\Re(\Omega)$

$$b_{\mathcal{D}}(\Pi_{FV}\boldsymbol{u}, \Pi_{\mathcal{C}}(\Pi_{FV}\boldsymbol{u})) \geq \Theta|u|^2_{\mathcal{T}}$$

Where the constant $\Theta$ is dependent on the mesh $\mathcal{D}$ but does not scale with $h$.

**Lemma 4.2:** Let $\mathcal{D}_n$ be a family of regular discretization triplets (in the sense that mesh parameters remain bounded) such that $h_n \to 0$, as $n \to \infty$. Furthermore, let $\theta_1$ and $\Theta$ be bounded independently of $n$. Then for all $n$, the solution $\boldsymbol{u}_n$ of Equation (4.1.6) exist and are unique, there exists $\widetilde{\boldsymbol{u}} \in (H_1(\Omega))^d$, and up to a subsequence (still denoted by $n$) $\Pi_{\mathcal{T}}\boldsymbol{u}_n \to \widetilde{\boldsymbol{u}}$ converges in $(L^q(\Omega))^d$, for $q \in [1, 2d/(d-2+\epsilon))$ as $h_n \to 0$. Finally, the cell-average finite volume gradient $\nabla_{\mathcal{D}}\boldsymbol{u}_n$ converges weakly to $\nabla\widetilde{\boldsymbol{u}}$ in $(L^2(\Omega))^{d^2}$.

**Theorem 4.3:** Consider the same case as in Lemma 4.2. Then the limit $\widetilde{\boldsymbol{u}} \in (H^1(\Omega))^d$ of the discrete mixed variational problem (4.1.6), and consequently the MPSA W-method, is the unique weak solution of Equation (4).

## 4.2 Comparison with earlier MPSA methods

The construction given in Section 3, as well as the analysis tools given in Section 4.1, allow us to give a qualitative comparison between the proposed MPSA methods. A quantitative comparison will follow in the results section.

First, we note that in analogy to the MPFA methods, so-called MPSA-U and MPSA-L methods can also be constructed [15]. However, these do not have a unique pointwise interpretation of stress, and as a consequence suffer from numerical locking for high Poisson ratios, and are not amenable to the analysis framework of Section 4.1. These methods are considered inferior in the context of elasticity, and will not be discussed further herein.

The MPSA-O method with strong symmetry does not lead to local problems with a unique solution given a single continuity point per face as in Equation (3.3.3) except on simplex grids with acute angles. To rectify this deficiency, a richer space of face variables needs to be introduced, where continuity is enforced weakly [15], [19], and the local linear system takes the form of a constrained minimization problem. The resulting method is termed the generalized MPSA-O method. This method is largely locking-free (except on PEBI grids), and shows good convergence properties. However, for simplex grids, the local coercivity condition does not in general hold, and the method performs relatively poorly on these grids.

The MPSA-W method proposed herein has a unique solution to the local problems under similar conditions to the MPFA-O, since they both have essentially equivalent locally coercive operators for the fine-scale problems. In practice, no examples have been found where the local problems do not have a unique solution, and as such, there is no need to introduce the weak displacement continuity used to obtain the generalized MPSA-O method. As a consequence, the method is significantly simpler, with much smaller local problems. We note that, should a case be found where the local problems are singular for the MPSA-W method, this can be handled by introducing multiple continuity points in (3.3.3) for this local problem, effectively constructing a generalized MPSA-W method. Our implementation can handle this option, however, we have little experience with the method.

An important observation in the analysis of MPFA and MPSA methods on simplex grids is that for the single continuity point chosen as O(1/3) (as defined in Section 3.3), the operators $\overline{\nabla}$ and $\widetilde{\nabla}$ coincide [24], [26]. For the MPFA-O method, this implies symmetry, and moreover the local coercivity condition is trivially satisfied. For the MPSA-W method, the same holds for homogeneous materials, and the method is therefore suitable for arbitrary simplex grids for homogeneous problems. None of the previous MPSA methods share this property.

## 5. Numerical experiments

In this section the performance of the MPSA-W method is compared to the generalized O-method. The quality measurements we consider are errors and convergence properties, computational cost and accuracy of the angular momentum on cells.

## 5.1 Convergence

To verify the performance of the MPSA-W method, we consider numerical examples which test the performance of the methods on homogeneous and heterogeneous problems in 2D and 3D. We consider Cartesian and simplex grids, and test both smooth grids and rough perturbations. To limit the influence of random perturbations, all variations of methods and material parameters are performed on the same perturbed grids.

### 5.1.1 Numerical setup

We first consider the unit square with a heterogeneity upper right corner, defined by the characteristic function

$$\chi_{2D}(x,y) = \begin{cases} 1, & \min(x,y) > 1/2 \\ 0, & otherwise \end{cases} \qquad (5.1.1)$$

The material parameters are then defined as $\mu = (1 - \chi_{2D}) + \kappa\chi_{2D}$, and $\lambda = \alpha\mu$; high values of $\alpha$ corresponds to a nearly incompressible material. For this medium, we consider the analytical solution

$$\boldsymbol{u}(x,y) = \begin{pmatrix} (x-0.5)^2(y-0.5)^2 \\ -\frac{2}{3}(x-0.5)(y-0.5)^3 \end{pmatrix} / ((1 - \chi_{2D}) + \kappa\chi_{2D}) \qquad (5.1.2)$$

The solution is depicted in Figure 2 together with the driving force. Note that the analytical solution is constructed so that the stress field is independent of the material heterogeneity.

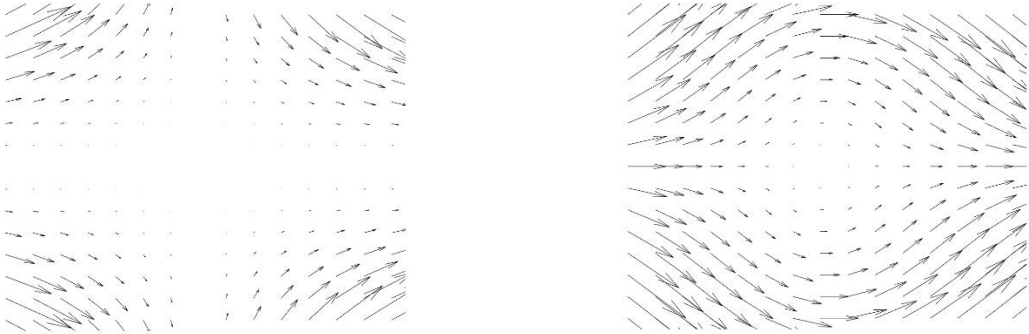

    a) Solution                                                   b) Driving force
Figure 2: Solution (left) and driving force (right) for the 2D problem in the homogeneous case ($\kappa = 1$).

We construct a similar three-dimensional solution on the unit cube. Again, we place a heterogeneity upper right corner, this time defined by the characteristic function

$$\chi_{3D}(x,y,z) = \begin{cases} 1, & \min(x,y,z) > 1/2 \\ 0, & otherwise \end{cases} \qquad (5.1.3)$$

The material parameters are then defined as $\mu = (1 - \chi_{3D}) + \kappa\chi_{3D}$, and $\lambda = \alpha\mu$. For this medium, we consider the analytical solution

$$\boldsymbol{u}(x,y,z) = \begin{pmatrix} (x-0.5)^2(y-0.5)^2(z-0.5)^2 \\ (x-0.5)^2(y-0.5)^2(z-0.5)^2 \\ -\frac{2}{3}((x-0.5)\ (y-0.5)^2 + (x-0.5)^2(y-0.5))(z-0.5)^3 \end{pmatrix} /((1-\chi_{3D}) + \kappa\chi_{3D})$$

(5.1.4)

This analytical solution is divergence free, with a stress field that is independent of the material heterogeneity.

To quantify the errors, we compute the discrete L2-norm of the displacement error as

$$\epsilon_D = \frac{\left(\sum_{K\in\mathcal{T}} m_K(\boldsymbol{u}(x_K) - \boldsymbol{u}_K)^2\right)^{\frac{1}{2}}}{\left(\sum_{K\in\mathcal{T}} m_K \boldsymbol{u}(x_K)^2\right)^{\frac{1}{2}}}$$

(5.1.5)

Similarly, the error in stresses are given by

$$\epsilon_\pi = \frac{\left(\sum_{\sigma\in\mathcal{F}} m_\sigma(\boldsymbol{\pi}(x_\sigma) - \boldsymbol{\pi}_\sigma)^2\right)^{\frac{1}{2}}}{\left(\sum_{\sigma\in\mathcal{F}} m_\sigma \boldsymbol{\pi}(x_\sigma)^2\right)^{\frac{1}{2}}}$$

(5.1.6)

### 5.1.2 Cartesian grids

Our first tests consider Cartesian grids. From Section 4, we expect both the generalized O-method and the new W-method to perform well in this case. Convergence results for 2D grids are reported in Figure 3. We observe that for regular grids, both methods are second order convergent in displacement, whereas the convergence rate for stress is reduced from roughly 1.5 on regular to first order on perturbed grids. We also observe that when a material heterogeneity is introduced by setting $\kappa = 10^6$, this has virtually no effect on the solution.

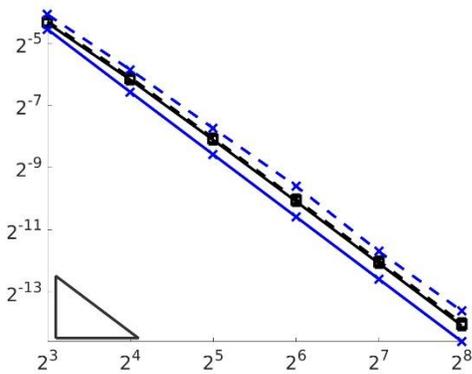
a) $\kappa = 1$, displacement

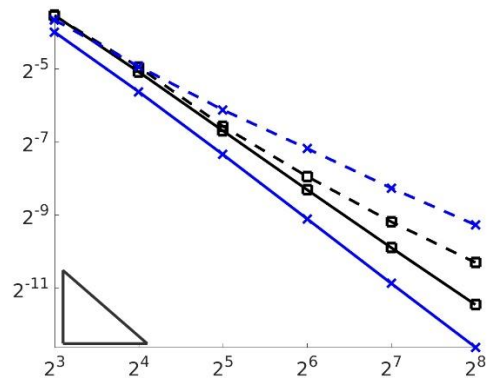
b) $\kappa = 1$, stress

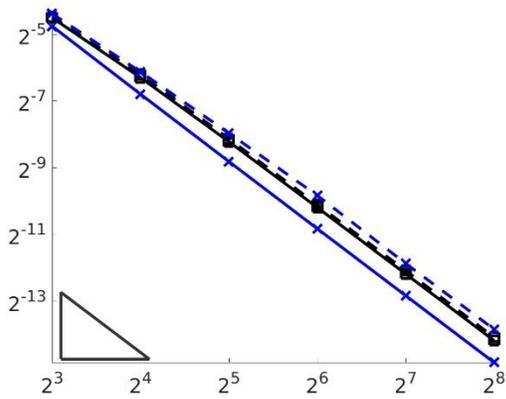
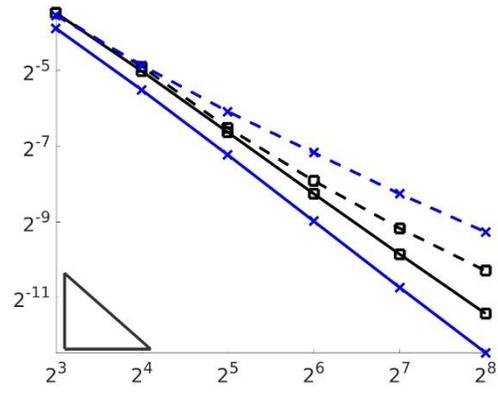

c) $\kappa = 10^6$, displacement

d) $\kappa = 10^6$, stress

Figure 3: Convergence on Cartesian grids in 2D for homogeneous and heterogeneous medium. Legend: Black lines with squares represent MPSA-W, blue lines with crosses the generalized method. Solid lines are regular grids, stippled lines perturbed grids. The axes show the square root of the number of degrees of freedom and the errors. The triangle indicates second order convergence.

Figure 4 shows the convergence behavior for the three dimensional test case. Again, all methods exhibit second order convergence in displacement, the generalized O-method actually has a slightly higher rate on unperturbed grids. The stresses have convergence rates of between 1.5 and 2 for unperturbed grids, and between 1 and 1.5 on perturbed grids. We also see that for the W-method, edge-wise and vertex-wise averaging renders very similar results, by close inspection the edge-wise approach is slightly better. Again, a material heterogeneity has only minor impact on the quality of the solutions.

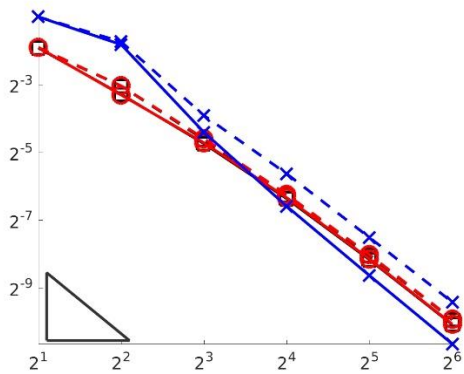
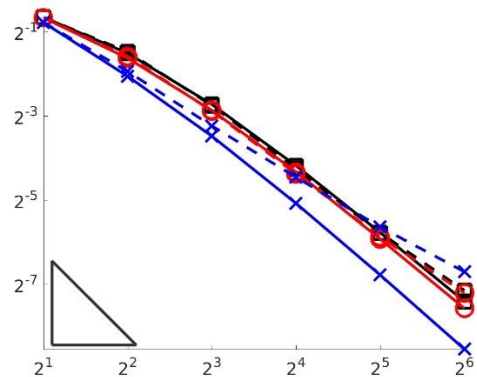

a) $\kappa = 1$, displacement

b) $\kappa = 1$, stress

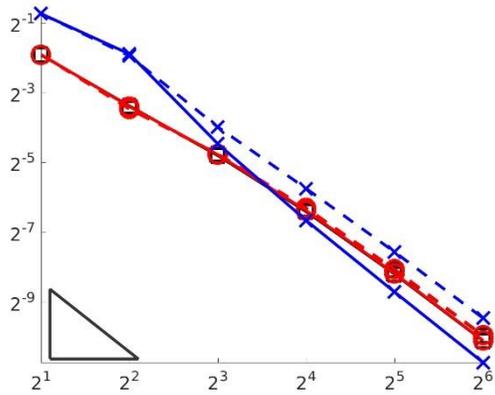
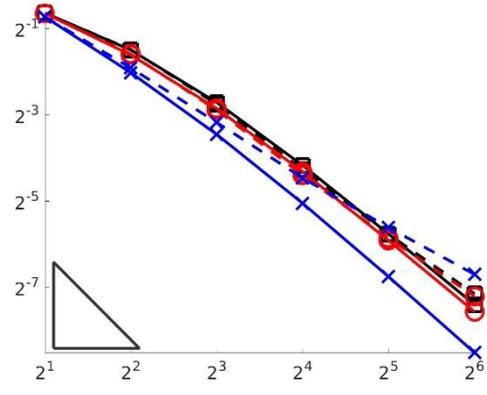

c) $\kappa = 10^6$, displacement   d) $\kappa = 10^6$, stress

Figure 4: Convergence on Cartesian grids in 3D for homogenous and heterogeneous media. Legend: Black lines with squares represent MPSA-W with vertex wise averaging, red lines with circles represent MPSA-W with edge wise averaging, blue lines with crosses the generalized method. Solid lines are regular grids, stippled lines perturbed grids. The axes show the cubic root of the number of degrees of freedom and the errors. The triangle indicates second order convergence.

### 5.1.3 Simplex grids

Next, we consider the same test examples, but now on simplex grids. As mentioned in Section 4, the generalized MPSA-O method does not produce a coercive local problem for these grids, and the method may deteriorate in this case. Conversely, the MPSA-W method is still coercive, and we expect better behavior of this method.

Figure 5 shows the convergence on triangular grids in 2D. The MPSA-W method yields second order convergence in the displacement, and around 1.5 for the stresses. This holds for both perturbed and unperturbed grids. Comparison with the corresponding test on Cartesian grids reveals that the errors are roughly the same for similar numbers of degrees of freedom. As for Cartesian grids, the convergence is robust with respect to material heterogeneities. On regular grids, the combination of the grid pattern and the loss of coercivity cause the generalized O-method to produce a resonance error, resulting in only first order convergence in displacement. Second order convergence is recovered when the grid is perturbed.

We also remark that the simplex grids employed in this test are created by dividing squares, rendering right angled triangles in the regular case. In additional tests on equilateral triangles (not shown here), we have observed that the generalized O-method may fail to converge at all. As noted in Section 4, the standard MPSA O-method with continuity point 1/3 fulfills the local coercivity condition in this case, and is therefore convergent. However, the standard O-method fails for local problems involving four cells meeting with right angles, creating a grid which is locally Cartesian. To apply the O-method on simplex grids in 2D, it is therefore necessary either to require acute angles only, or to switch between the standard and generalized method depending on the local geometry. The MPSA-W method offers an attractive alternative which avoids these issues altogether.

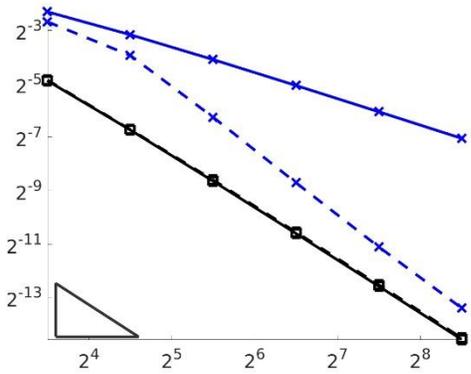

a)  $\kappa = 1$, displacement

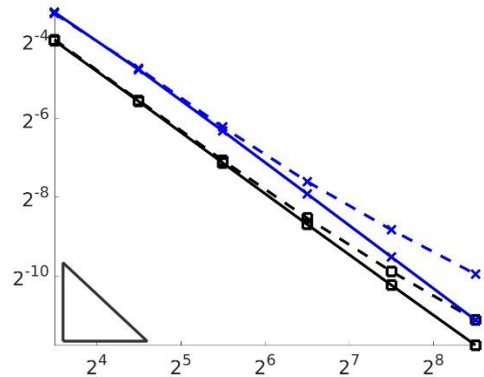

b)  $\kappa = 1$, stress

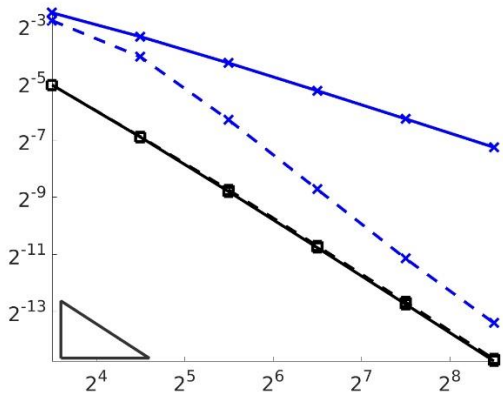

c)  $\kappa = 10^6$, displacement

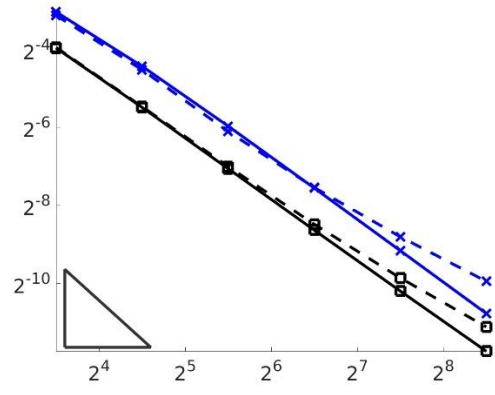

d)  $\kappa = 10^6$, stress

Figure 5: Convergence on simplex grids in 2D for homogeneous and heterogeneous medium. Legend: Black lines with squares represent MPSA-W, blue lines with crosses the generalized method. Solid lines are regular grids, stippled lines perturbed grids. The axes show the square root of the number of degrees of freedom and the errors. The triangle indicates second order convergence.

For tetrahedral grids, the W-methods again perform well, with second order convergence on displacements, and between first and second for stresses, see Figure 6. As for the Cartesian grids, edge-wise averaging consistently produces somewhat smaller errors, however the difference is minor. The generalized O-method only slowly approaches its convergence region; on perturbed grids the stresses do not converge for the grid refinements considered here. We observe that in this case, the material heterogeneity actually improves the convergence behavior.

While the generalized O-method is used as a benchmark herein, a comparison between the standard O-method and lowest order finite elements (FEM) was undertaken in [15]. For all cases considered therein, the O-method was at least as good as the FEM method both in terms of displacements and stresses. We have not carried out a similar comparison here, but the experiments in this paper indicate that the MPSA-W method will also compare favorably with the FEM.

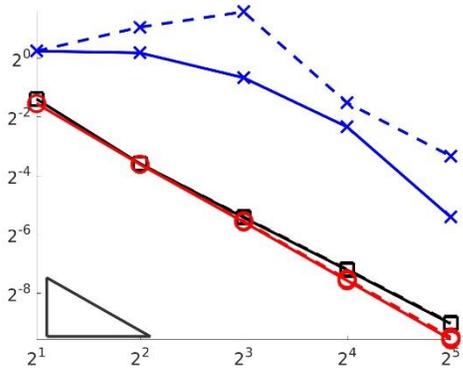
e)  $\kappa = 1$, displacement

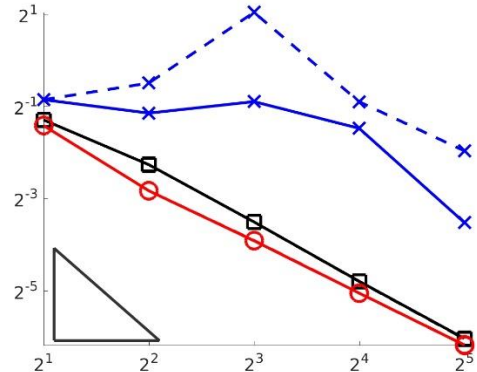
f)  $\kappa = 1$, stress

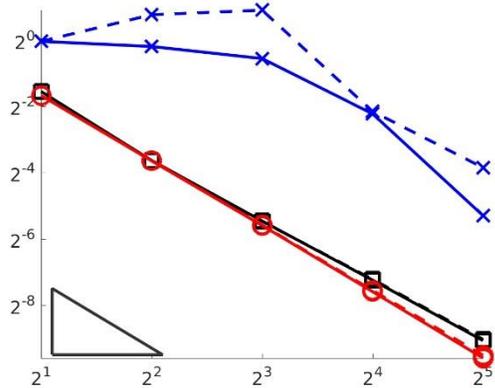
g)  $\kappa = 10^6$, displacement

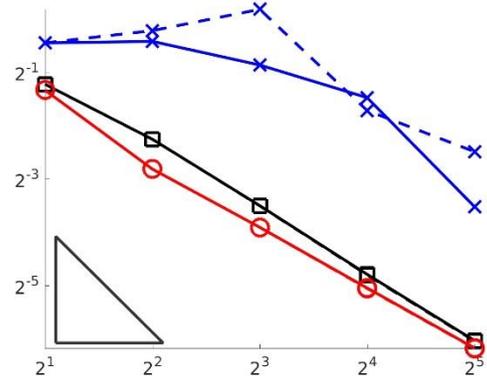
h)  $\kappa = 10^6$, stress

Figure 6: Convergence on simplex grids in 3D for homogeneous and heterogeneous medium. Legend: Black lines with squares represent MPSA-W with vertex wise averaging, red lines with circles represent MPSA-W with edge wise averaging, blue lines with crosses the generalized method. Solid lines are unperturbed grids, stippled lines perturbed grids. The axes show cubic root of the number of degrees of freedom and the errors. The triangle indicates second order convergence.

### 5.1.4 Robustness to volume locking

The MPSA methods have been shown to be free from volume locking for simplex grids [19], whereas for Cartesian girds, there is no proof that the methods will not exhibit locking. To test the methods in the incompressible limit, we ran both the 2D and 3D simulations with the ratio $\alpha = \lambda/\mu = \{10^2, 10^3, 10^4\}$, corresponding to Poisson ratios $\nu = \lambda/(2(\lambda + \mu))$ of 0.495, 0.4995 and 0.49995, respectively. Only results from the 3D simulations are reported here, the 2D case showed a similar, although on average somewhat better behavior. For simplicity we consider a homogeneous domain, parameter heterogeneities did not alter the convergence properties.

Figure 7 depicts the convergence ratios on a Cartesian grid. We observe that for unperturbed grids, both displacement and stresses exhibit steady convergence, with rates around 2 for displacement and around 1.5 for stresses. The generalized method yields a somewhat smaller error than the MPSA-W formulations, in particular for stresses. As the grid is perturbed, the displacement is still convergent, although with an irregular behavior for the generalized method for the highest Poisson ratio. The convergence rate for stresses is more sensitive to perturbations, in particular for the generalized method, and for extreme Poisson ratios, it is questionable whether the method can be said to converge at all. We also note that for these cases, edge wise and vertex wise averaging of the stresses yields virtually indistinguishable results.

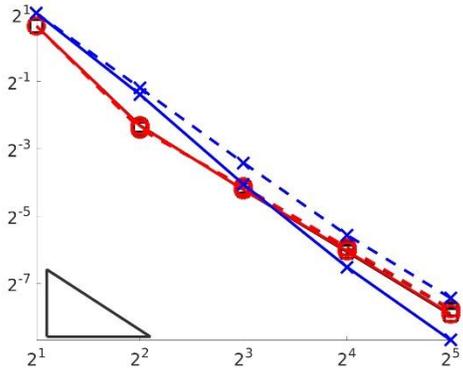

a) $\nu = 0.495$, displacement

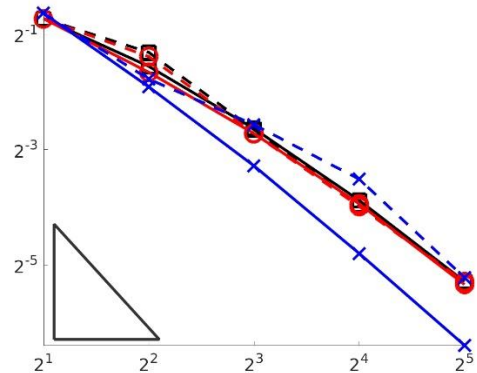

b) $\nu = 0.495$, stress

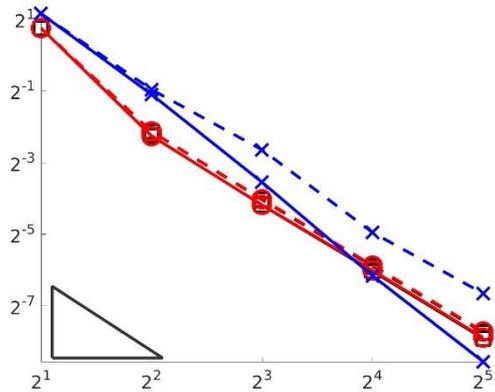

c) $\nu = 0.4995$, displacement

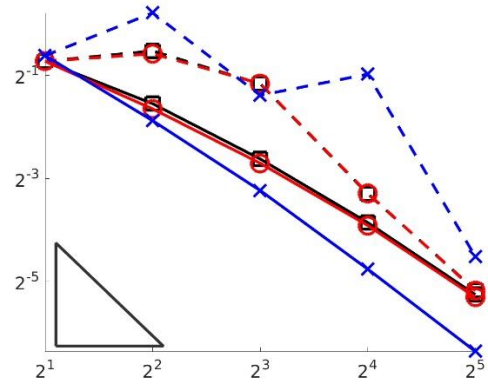

d) $\nu = 0.4995$, stress

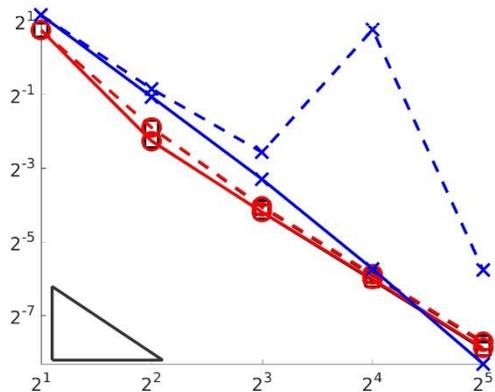

e) $\nu = 0.49995$, displacement

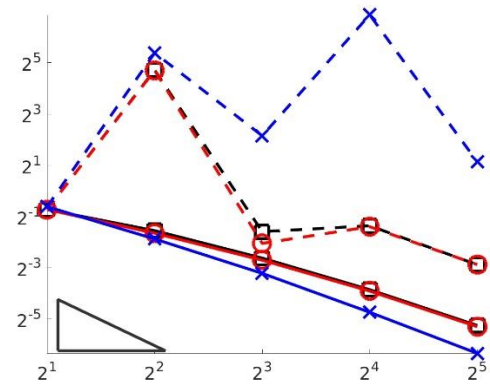

f) $\nu = 0.49995$, stress

Figure 7: Homogeneous domain in 3D with Cartesian grid and increasing Poisson ratio. Legend: Black lines represent MPSA-W with vertex wise averaging, dark gray lines the generalized method, light gray

lines MPSA-W with edgewise averaging. Solid lines are regular grids, stippled lines perturbed grids. The axes show square root of the number of degrees of freedom and the errors. The triangle indicates second order convergence.

Next, we consider the performance on tetrahedral grids, shown in Figure 8. As expected, the convergence behavior is better than for the Cartesian grid, both for unperturbed and perturbed grids, and the methods shows second order convergence for displacement and first order for the stresses, independent of the Poisson ratio. The only exception is for stresses on perturbed grids with $\nu = 0.49995$, where the MPSA-W methods show irregular convergence, and it is not clear whether the generalized method enters the region of convergence at all. We remark that for high Poisson, the iterative linear solvers available to us provided solutions that were inaccurate and polluted the convergence rate of the discretization schemes, even for strict residual tolerances. This forced us to rely upon direct solvers for these cases (as implemented by the Matlab backslash operator), which limited how many refinement steps were feasible.

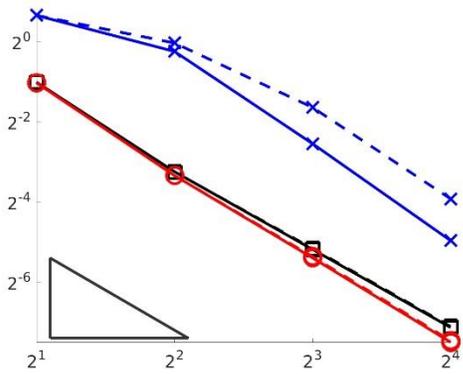

a) $\nu = 0.495$, displacement

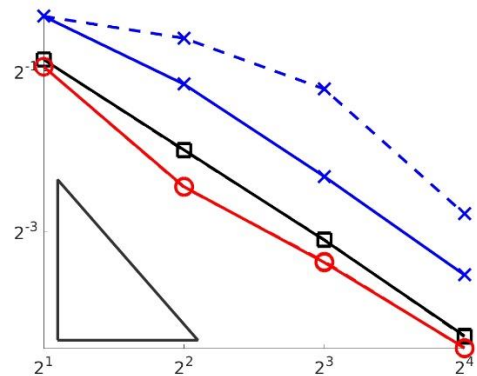

b) $\nu = 0.495$, stress

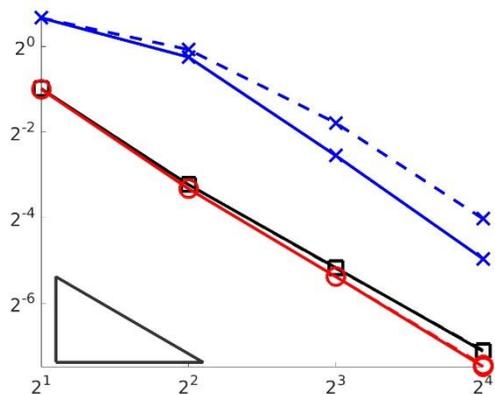

c) $\nu = 0.4995$, displacement

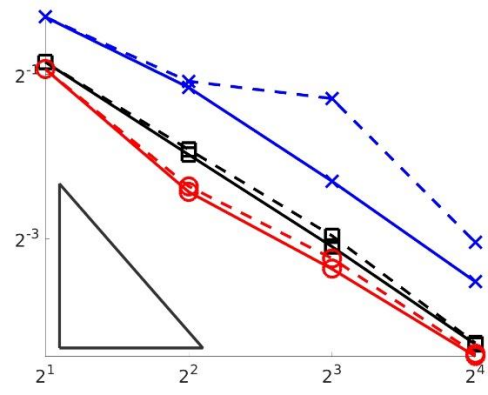

d) $\nu = 0.4995$, stress

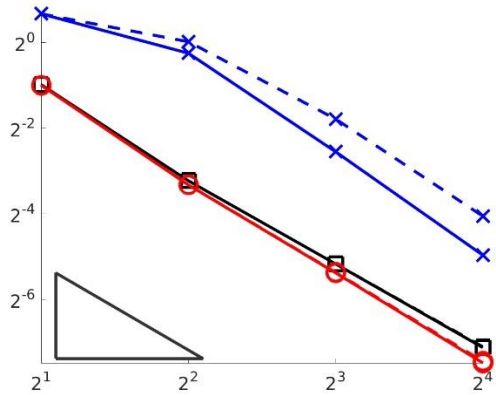 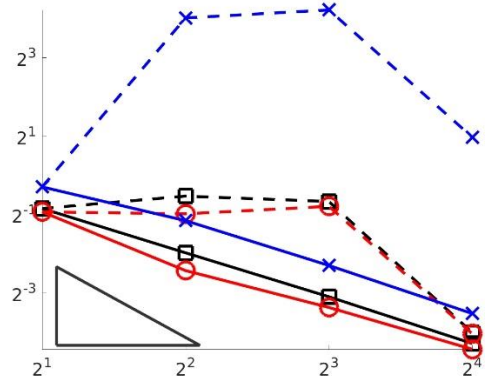

e)  $\nu = 0.49995$, displacement     f)  $\nu = 0.49995$, stress

Figure 8: Homogeneous domain in 3D with tetrahedral grid and increasing Poisson ratio. Legend: Black lines represent MPSA-W with vertex wise averaging, dark gray lines the generalized method, light gray lines MPSA-W with edgewise averaging. Solid lines are unperturbed grids, stippled lines perturbed grids. The axes show the cubic root of the number of degrees of freedom and the errors. The triangle indicates second order convergence.

## 5.2 Angular momentum

Since the W-method imposes the symmetry weakly, it is of interest to consider whether the method produces an asymmetric stress tensor. While the MPSA methods only provide stresses on the faces, it is possible to reconstruct a point-wise stress tensor by solving Neumann problems on each cell, with the discrete stresses as boundary conditions. This post-processed stress tensor will be locally symmetric both for the MPSA-O and MPSA-W methods.

A related question is whether the weak symmetry impacts the accuracy of the rotational components of the stress. As a final test of the method, we therefore consider the angular momentum of the surface forces on each cell. Denote by $\boldsymbol{T}^t_{K,\sigma}$ the tangent component of the stresses on surface $\sigma$, as seen from cell $K$. The angular momentum of $\boldsymbol{T}^t_{K,\sigma}$ around the cell center is given by $\omega_{K,\sigma} = \boldsymbol{T}^t_{K,\sigma} \times (\boldsymbol{x}_\sigma - \boldsymbol{x}_K)$, and the total angular momentum on $K$ is $\omega_K = \sum_{\sigma \in \mathcal{F}_K} \omega_{K,\sigma}$. We then measure the quality of the approximated angular momentum by the mean of the absolute value over all cells in the grid.

Two metrics were considered for the angular momentum: Frist, due to symmetry of the stress tensor, under grid refinement, the angular momentum should converge with a rate of more than $h^D$ ($h$ being the grid diameter), and this was verified for all cases considered. This is the analogue to the Cauchy stress theorem. Second, we measured the difference in angular momentum between the numerical and analytical surface stresses. Figure 9 shows the convergence of this quantify for the 2D and 3D test problems with a homogenous domain ($\kappa = 1$). The angular momentum for both the O- and W-method converge with a rate of above 3, and despite the weakly imposed symmetry, the W-method is as accurate as the O-method.

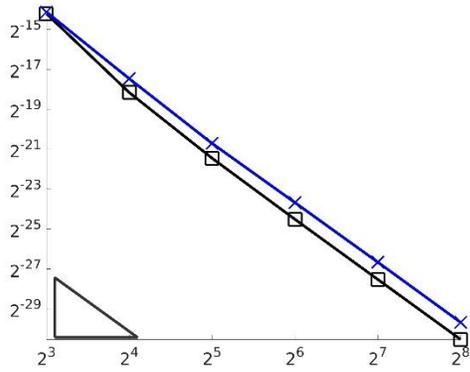
a) Cartesian – 2D

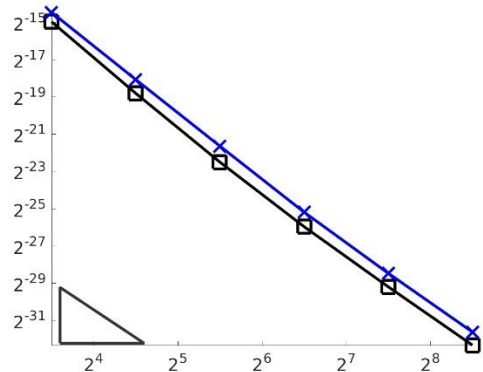
b) Triangular – 2D

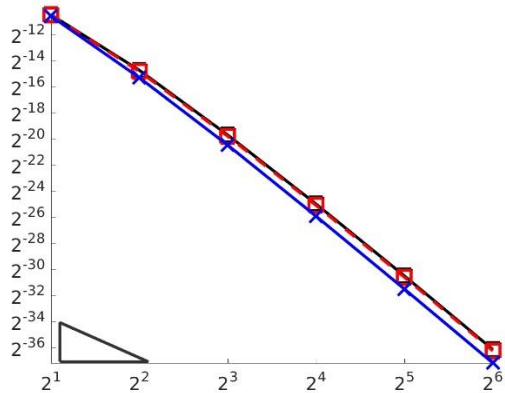
c) Cartesian – 3D

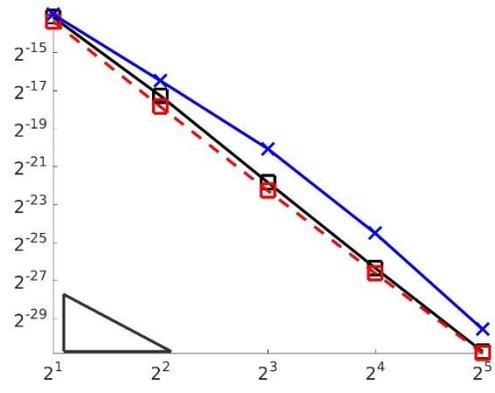
d) Tetrahedral – 3D

Figure 9: Difference between analytical and numerical angular momentum around cell centers, taken as mean value over all cells. Legend: W-method node-wise in black squares, W-Method edge-wise in red squares (only 3D), black crosses: generalized method. The axes show square (2D) or cubic (3D) root of the number of degrees of freedom and the errors. Triangles indicate third order convergence.

## 5.3 Computational cost

The main computational cost in the discretization is related to the inversion of the linear system (3.3.4). Here the MPSA-W methods gives a significant computational saving compared to the generalized approach: First, the number of displacement continuity points is significantly reduced. This is particularly important in 3D, where the number of displacement continuity Equations is reduced from 144 to 36 for Cartesian grids, and on average from about 360 to 90 for simplex grids considered herein. Second, in the weakly symmetric formulation, the linear system is no longer overdetermined, and it can be solved as a direct problem, rather than a constrained optimization problem.

The MPSA methods considered herein all lead to the same cell stencils, and can be expected to have a similar cost in solving the global linear system. While a systematic investigation is beyond the scope of this paper, it is our experience that when iterative solvers are employed, specifically GMRES with ILU and AMG preconditioners, the number of iterations are roughly the same for all methods. The exception of the O-method on simplex grids, where the 2-4 times more iterations may be needed. Again, we suspect this is related to a lack of local coercivity.

## 6. Concluding remarks

In this paper, we have introduced a weakly symmetric finite volume method for elasticity, termed the W-method. The method uses cell-centered variables, and it is particularly well suited for simulations of coupled flow and deformation in porous media. The MPSA-W method differs from previous MPSA methods in that it is robust on both simplex and Cartesian grids, and it has a lower computational cost.

Convergence of the W-method was shown by casting the method in a discrete variational form previously established for the MPSA-O method. The convergence criteria, which depend only on the grid and the constitutive law, can be verified at the time of discretization. While both the W- and O-methods are convergent on Cartesian grids, only the W-method has this property also on simplex grids, whereas the O-method loses local coercivity.

We present an extensive suite of numerical experiments, including the first comprehensive convergence tests for the MPSA methods in 3D. Both the W- and O-method performed well on Cartesian grids, but only the W-method showed similar behavior on simplex grids. While the weak symmetry in the W-method in 3D can be imposed by vertex-wise or edge-wise averaging, the numerical experiments further indicated that the extra accuracy obtained by the edge-wise approach do not justify the additional data structure needed, and thus we prefer vertex-wise averaging.

The results reported herein allow us to conclude that the new MPSA-W method represents the most simple, efficient and robust version of the MPSA methods currently available, and forms the benchmark for future developments of the methods.